\documentclass[11pt,psamsfonts]{amsart}
\usepackage{amsmath}
\usepackage{amsthm}
\usepackage{amssymb}
\usepackage{amscd}
\usepackage{amsfonts}
\usepackage{wrapfig}
\usepackage{amsbsy}
\usepackage{epsfig,afterpage}
\usepackage[dvips]{psfrag}
\usepackage[all]{xy}
\usepackage{pstcol}
\newcommand{\R}{\ensuremath{\mathbb{R}}}

\newcommand{\N}{\ensuremath{\mathbb{N}}}

\newtheorem {definition}  {Definition}
\newtheorem {remark}  {Remark}
\newtheorem {example} {Example}

\definecolor{red}{rgb}{1.,0.,0.}
\definecolor{blue}{rgb}{0.,0.,1.}
\definecolor{pink}{rgb}{1.,0.75,0.8}

\begin{document}

\title[Fold--Saddle Bifurcation in Non$-$Smooth Systems] {Fold$-$Saddle Bifurcation in Non--Smooth Vector Fields on the Plane.}

\author[C.A. Buzzi, T. de Carvalho and M.A. Teixeira]
{Claudio A. Buzzi$^1$, Tiago de Carvalho$^1$ and\\ Marco A.
Teixeira$^2$}

\address{$^1$ IBILCE--UNESP, CEP 15054--000,
S. J. Rio Preto, S\~ao Paulo, Brazil}

\address{$^2$ IMECC--UNICAMP, CEP 13081--970, Campinas,
S\~ao Paulo, Brazil}

\email{buzzi@ibilce.unesp.br}

\email{tiago@ibilce.unesp.br}

\email{teixeira@ime.unicamp.br}

\subjclass{Primary 34A36, 37G10, 37G05}

\keywords{Fold$-$Saddle singularity, canard, limit cycle,
bifurcation, non$-$smooth vector field.}
\date{}
\dedicatory{} \maketitle


\begin{abstract}

 This paper presents  results concerning bifurcations of $2D$
piecewise$-$smooth dynamical systems governed by vector fields.
Generic three$-$parameter families of a class of Non$-$Smooth Vector
Fields are studied and its bifurcation diagrams are exhibited. Our
main results describe the unfolding of the so called $Fold-Saddle$
singularity.

\end{abstract}

\section{Introduction}

The general purpose of this article is to present some aspects of
the geometric and qualitative theory of a class of planar
non$-$smooth systems. Our main concern is to discuss the behavior of
such systems around typical singularities that appear generically in
three$-$parameter families. We mention that certain phenomena in
control systems, impact in mechanical systems and nonlinear
oscillations are the main sources of motivation of our study
concerning the dynamics of those systems that emerge from
differential equations with discontinuous right$-$hand sides.

The codimension zero and codimension one singularities were
discussed in \cite{K} and \cite{Kuznetsov} respectively. In
\cite{Marcel} codimension two singularities were studied. 
The specific topic addressed in this paper is
the complete characterization of the \textit{Fold$-$Saddle
bifurcation diagram}. Those papers give the necessary
basis for the development of our approach.\\

Let $K\subseteq \R ^{2}$ be a compact set and $\Sigma \subseteq K$
given by $\Sigma =f^{-1}(0),$ where $f:K\rightarrow \R$ is a smooth
function having $0\in \R$ as a regular value (i.e. $\nabla f(p)\neq
0$, for any $p\in f^{-1}({0}))$ such that $\partial K \cap \Sigma =
\emptyset$ or $\partial K \pitchfork \Sigma$. Clearly $\Sigma$ is
the separating boundary of the regions $\Sigma_+=\{q\in K | f(q)
\geq 0\}$ and $\Sigma_-=\{q \in K | f(q)\leq 0\}$. We can assume
that $\Sigma$ is represented, locally
around a point $q=(x,y)$, by the function $f(x,y)=y.$\\

Designate by $\chi^r$ the space of $C^r$ vector fields on $K$
endowed with the $C^r-$topology with $r\geq 1$ or $r=\infty,$ large
enough for our purposes. Call \textbf{$\Omega^r=\Omega^r(K,f)$} the
space of vector fields $Z: K \setminus\Sigma \rightarrow \R ^{2}$
such that
$$
 Z(x,y)=\left\{\begin{array}{l} X(x,y),\quad $for$ \quad (x,y) \in
\Sigma_+,\\ Y(x,y),\quad $for$ \quad (x,y) \in \Sigma_-,
\end{array}\right.
$$
where $X=(f_1,g_1)$, $Y = (f_2,g_2)$ are in $\chi^r.$ We write
$Z=(X,Y),$ which we will accept to be multivalued in points of
$\Sigma.$ The trajectories of $Z$ are solutions of  ${\dot q}=Z(q),$
which has, in general, discontinuous right$-$hand side. The basic
results of differential equations, in this context, were stated by
Filippov in
\cite{Fi}. Related theories can be found in \cite{K, ST, T}.\\

In what follows we will use the notation
\[X.f(p)=\left\langle \nabla f(p), X(p)\right\rangle \quad \mbox{ and } \quad Y.f(p)=\left\langle \nabla f(p), Y(p)\right\rangle. \]


\subsection{Setting the problem}

Let $X_0$ be a smooth vector field defined in $\Sigma_{+}$. We say
that a point $p_0 \in \Sigma$ is a \textit{$\Sigma-$fold point} of
$X_0$ if $X_0.f(p_0)=0$ but $X_{0}^{2}.f(p_0)\neq0.$ Moreover,
$p_0\in\Sigma$ is a \textit{visible} (respectively {\it invisible})
$\Sigma-$fold point of $X_0$
 if $X_0.f(p_0)=0$ and $X_{0}^{2}.f(p_0)> 0$
(resp. $X_{0}^{2}.f(p_0)< 0$). In this universe,
$\Gamma^{X_0}_{\Sigma}$, a $\Sigma-$fold point has codimension zero.
Since $f(x,y)=y$ we derive the following generic normal forms
$X_0(x,y)= (\alpha_1, \beta_1 x)$ with $\alpha_1=\pm1$ and $\beta_1=
\pm1$.

Let $Y_0$ be a smooth vector field defined in $\Sigma_{-}$. Assume
that $Y_0$ has a hyperbolic saddle point $S_{Y_0}$ on $\Sigma$ and
that  the eigenspaces of $DY_0(S_{Y_0})$ are transverse to $\Sigma$
at $S_{Y_0}$. In this universe, $\Gamma^{Y_0}_{\Sigma}$, a saddle
point $S_{Y_{0}}$ has codimension one. Since $f(x,y)=y$ we derive
the following generic normal forms $Y_0(x,y)= (\alpha_2 y, \alpha_2
x)$ with $\alpha_2=\pm1$ and its generic unfolding $Y_{\beta}= (
\alpha_2 (y + \beta) ,  \alpha_2  x)$ where $\beta \in \R$. Let $U$
be a small neighborhood of $Y_0$ in $\Gamma^{Y_0}_\Sigma.$ Then:

\noindent (a) There exists a smooth function $L: U\rightarrow \R$,
such that $DL_{Y_0}$ is surjective.

\noindent (b) The correspondence $Y\rightarrow S_Y$ is smooth, where
$S_Y$ is a saddle point of $Y$.

\noindent (c) If $L(Y) >0$ then $S_Y \in \Sigma_-$.

\noindent (d) If $L(Y) =0$ then $S_Y \in \Sigma$.

\noindent (e) If $L(Y) <0$ then $S_Y  \in \Sigma_+$.

In this paper we are concerned with the bifurcation diagram of
systems  $Z_0=(X_0,Y_0)$ in $\Omega^r$ such that $p_0=S_{Y_0} \in
\Sigma$. This singularity will be called $\mathbf{Fold - Saddle}$
singularity (see Figures \ref{fig fold-sela} and \ref{fig
fold-visivel-sela}).

\begin{figure}[!h]
\begin{minipage}[b]{0.49\linewidth}
\begin{center}
 \epsfxsize=3.5cm \epsfbox{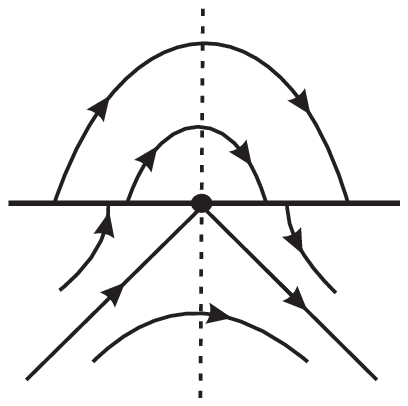}
\caption{\small{(Invisible) Fold-Saddle Singularity.}} \label{fig
fold-sela}
\end{center}
\end{minipage} \hfill
\begin{minipage}[b]{0.49\linewidth}
\begin{center}
 \epsfxsize=3.5cm \epsfbox{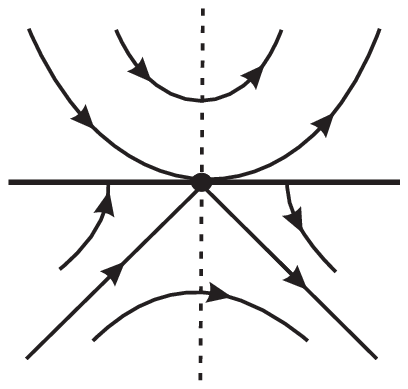}
\caption{\small{(Visible) Fold-Saddle Singularity.}} \label{fig
fold-visivel-sela}
\end{center}
\end{minipage}\end{figure}

We depart from $Z^{i}_0, Z^{v}_0 \in \Omega^r$ written in the
following forms:
\begin{equation}\label{eq fold-sela inicio}
Z^{i}_0 = \left\{
      \begin{array}{ll}
        X^{i}_0 = \left(
              \begin{array}{c}
                1 \\
               -x
\end{array}
      \right)
 & \hbox{if $y \geq 0$,} \\
        Y_0 =  \left(
              \begin{array}{c}
                -y \\
               -x
\end{array}
      \right)& \hbox{if $y \leq 0$, and}
      \end{array}
    \right.
\end{equation}
\begin{equation}\label{eq fold-visivel-sela inicio}
Z^{v}_0 = \left\{
      \begin{array}{ll}
        X^{v}_0 = \left(
              \begin{array}{c}
                1 \\
               x
\end{array}
      \right)
 & \hbox{if $y \geq 0$,} \\
        Y_0 =  \left(
              \begin{array}{c}
                -y \\
               -x
\end{array}
      \right)& \hbox{if $y \leq 0$.}
      \end{array}
    \right.
\end{equation}
Note that $X^{i}_0$ presents an invisible $\Sigma-$fold point on its
phase portrait and $X^{v}_0$ presents a visible one. Following the
techniques developed in \cite{T1}, we are able to prove that there
exists a smooth mapping
 $F_{\tau}:\Omega^r,Z^{\tau}_0 \rightarrow \R^{3},0$ where $\tau=i$ or $v$  
  such that:

1- $(DF_\tau)_{Z^{\tau}_0}$ is surjective (So
$M_{\tau}=(F_{\tau})^{-1}(0)$ is locally, around $Z^{\tau}_0,$ an
imbedded differentiable manifold).

2- Each $Z \in U_{\tau}$, with $F_{\tau}(Z)=0$ and $U_{\tau}$ a
small neighborhood of $Z^{\tau}_0$ in $\Omega^r$ is $C^0-$equivalent
to $Z^{\tau}_0.$

The main question is to exhibit the bifurcation diagram of
$Z^{\tau}_0$. So, we have to consider generic imbeddings
$\sigma_{\tau}: \R^{3},0 \rightarrow \Omega^r,Z^{\tau}_0$
($3-$parameter families). They are transversal imbeddings to
$M_{\tau}$ at $Z^{\tau}_0.$

\begin{figure}[!h]
\begin{minipage}[b]{0.31\linewidth}
\begin{center}
\epsfxsize=3cm  \epsfbox{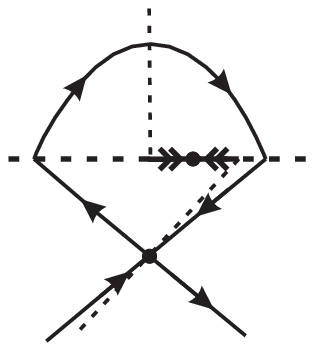}
\caption{\small{}}\label{fig perto estavel}
\end{center}
\end{minipage} \hfill
\begin{minipage}[b]{0.31\linewidth}
\begin{center}
 \epsfxsize=3cm \epsfbox{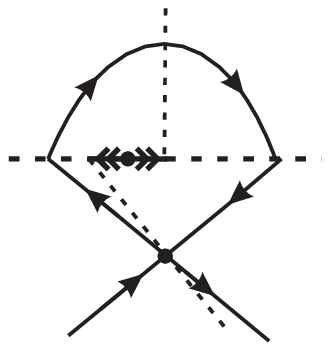}
\caption{\small{}} \label{fig perto instavel}
\end{center}
\end{minipage} \hfill
\begin{minipage}[b]{0.35\linewidth}
\begin{center}
 \epsfxsize=3cm \epsfbox{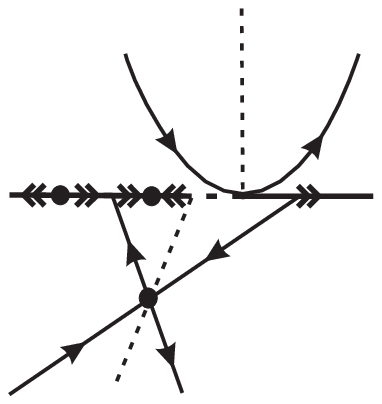}
\caption{\small{}} \label{fig sela fold generica}
\end{center}
\end{minipage}
\end{figure}

Consider $Z^{\tau}_0=(X^{\tau}_0,Y_0) \in U_{\tau}$. Roughly
speaking, we derive that:

I- There is a canonical imbedding $F_{0}^{\tau}:\R^2,0 \rightarrow
\chi^r,Z^{\tau}_0$ such that
$F_{0}^{\tau}(\lambda,\beta)=Z^{\tau}_{\lambda,\beta}$ expressed by:
\begin{equation}\label{eq fold-sela sem variar}
Z^{\tau}_{\lambda, \beta} = \left\{
      \begin{array}{ll}
        X^{\tau}_{\lambda} = \left(
              \begin{array}{c}
                   1 \\
                \alpha_1(\tau)(x - \lambda)
\end{array}
      \right)
 & \hbox{if $y \geq 0$,} \\
         Y_{\beta} = \left(
              \begin{array}{c}
               -(y + \beta) \\
               - x
\end{array}
      \right)& \hbox{if $y \leq 0,$}
      \end{array}
    \right.
\end{equation}
where   $\lambda$, $\beta \in (-1,1)$, $\alpha_1(i)=-1$ and
$\alpha_1(v)=1$. Moreover, its bifurcation diagram of
$Z^{\tau}_{\lambda, \beta}$ is exhibited (see Figures \ref{fig
diagrama bif teo 1} and \ref{fig diagrama bif teo 4 5 e 6}). We
observe that there are some typical topological types nearby
$Z^{\tau}_0$ that do not appear in the bifurcation diagram of
$Z^{\tau}_{\lambda,\beta}$. For example, when $\tau=i$ the
configurations in Figures \ref{fig perto estavel} and \ref{fig perto
instavel} are excluded and when $\tau=v$ the configuration in Figure
\ref{fig sela fold generica} also is excluded.

II- We add an auxiliary parameter $\mu$  in the following way:
\begin{equation}\label{eq fold-sela com parametros}
\overline{Z}^{\tau}_{\lambda , \mu, \beta} = \left\{
      \begin{array}{ll}
        X_{\lambda} = \left(
              \begin{array}{c}
                  1 \\
                \alpha_1(\tau)(x - \lambda)
\end{array}
      \right)
 & \hbox{if $y \geq 0$,} \\
         Y_{\mu,\beta} = \left(
              \begin{array}{c}
                \frac{\mu}{2} x + \frac{(\mu - 2)}{2}(y + \beta) \\
                \frac{(\mu - 2)}{2} x + \frac{\mu}{2}(y + \beta)
\end{array}
      \right)& \hbox{if $y \leq 0,$}
      \end{array}
    \right.
\end{equation}
where   $\lambda$, $\beta \in (-1,1)$, $\alpha_1(i)=-1$,
$\alpha_1(v)=1$ and $\mu \in (- \varepsilon_0, \varepsilon_0)$ with
the real number $\varepsilon_0 > 0$ being sufficiently small. By
means of this late unfolding its bifurcation diagram cover all
topological types near $\overline{Z}^{\tau}_{0,0,0}$.

In this universe, $\Gamma^{Z^{\tau}_0}$, a Fold$-$Saddle singularity
has codimension three. Since $f(x,y)=y$ we derive the generic normal
forms $\overline{Z}^{\tau}_{\lambda , \mu_0, \beta}$ with $\mu_0=\pm
\varepsilon_0/2 $ and its generic unfolding
$\overline{Z}^{\tau}_{\lambda , \mu,
\beta}$ given by \eqref{eq fold-sela com parametros}
. Therefore, there is a codimension three bifurcation
(global) branch terminating at $Z^{\tau}_0$. In fact, note that we
can obtain Equation \eqref{eq fold-sela inicio} (respectively
\eqref{eq fold-visivel-sela inicio}) from Equation \eqref{eq
fold-sela com parametros}  taking $\tau=i$ (respectively $\tau=v$),
$\lambda=0$, $\mu=0$ and $\beta=0$ .


%
%
%

%

Of course, we can take another generic normal form of one or both
vector fields $X_0$ and $Y_0$. In this paper we consider just the
cases described in Equations \eqref{eq fold-sela inicio} and
\eqref{eq fold-visivel-sela inicio}. For the other cases a similar
approach can be done.

%

It is worth mentioning that  we detect branches of ``\emph{canard
cycles}"  in the bifurcation diagram of
$\overline{Z}^{i}_{\lambda,\mu,\beta}$. Recall that, a canard cycle
is a closed path composed by pieces of orbits of $X$, $Y$ and
$Z^\Sigma$ (see Figures \ref{fig canard I}, \ref{fig canard II} and
\ref{fig canard}). In Section 2  a precise definition will be given.

\begin{example} Equations
\eqref{eq fold-sela inicio} and \eqref{eq fold-visivel-sela inicio}
appear in problems related to \textit{Control Theory}, more
specifically, in \textit{Relay Systems}. In fact, consider the
function $\varphi: \R \rightarrow \R$ given by $$\varphi(y)= \left\{
                                                               \begin{array}{rl}
                                                                 -1, & \hbox{for $y\geq0$,} \\
                                                                 -y, & \hbox{for $y\leq0$,}
                                                               \end{array}
                                                             \right.
$$
and $u(y)=-\varphi(y) sign (y)$. So \eqref{eq fold-sela inicio} and
\eqref{eq fold-visivel-sela inicio} are represented by
$\widetilde{Z}^{\tau}_{0}(x,y) =$ $(u(y),\overline{\alpha}(\tau) x)$
where $\tau=i$ or $v$, $\overline{\alpha}(i)=-1$ and
$\overline{\alpha}(v)=sign (y)$.
\end{example}

\subsection{Statement of the Main Results} Our results are now
stated. Theorems 1, 2 and 3  are intermediate steps towards Theorem
A and Theorems 4, 5 and 6 are intermediate steps towards Theorem
B.\vspace{.5cm}


\noindent {\bf Theorem 1.} {\it Take $\tau=i$ in Equation \eqref{eq
fold-sela sem variar} or equivalently, take $\tau=i$ and $\mu=0$ in
Equation \eqref{eq fold-sela com parametros}. The
$(\lambda,\beta)-$plane contains essentially $17$ distinct typical
configurations representing $5$ distinct topological behaviors on
its bifurcation diagram (see Figure \ref{fig diagrama bif teo 1}).}

\vspace{.5cm}

It is easy to see that the cases covered by Theorem 1 do not
represent the full unfolding of the (Invisible) Fold$-$Saddle singularity. 
Because of this, 
the next two theorems are necessary. Each one of them describes a
distinct generic codimension two
singularity.\\

\noindent {\bf Theorem 2.} {\it Take $\tau=i$ and $0 <\mu <
\varepsilon_0$ in Equation \eqref{eq fold-sela com parametros}. The
$(\lambda,\beta)-$plane contains essentially $19$ distinct typical
configurations representing $7$ distinct topological behaviors on
its bifurcation diagram (see Figure \ref{fig diagrama bif teo 2}).}

\vspace{.5cm}

\noindent {\bf Theorem 3.} {\it Take $\tau=i$ and
$-\varepsilon_0<\mu < 0$ in Equation \eqref{eq fold-sela com
parametros}. The $(\lambda,\beta)-$plane contains essentially $19$
distinct typical configurations representing $7$ distinct
topological behaviors on its bifurcation diagram (see Figure
\ref{fig diagrama bif teo 3}).}

\vspace{.5cm}

\noindent {\bf Theorem 4.} {\it Take $\tau=v$ in Equation \eqref{eq
fold-sela sem variar} or equivalently, take $\tau=v$ and $\mu=0$ in
Equation \eqref{eq fold-sela com parametros}. The
$(\lambda,\beta)-$plane contains essentially $13$ distinct typical
configurations representing $7$ distinct topological behaviors on
its bifurcation diagram (see Figure \ref{fig diagrama bif teo 4 5 e
6}).}

\vspace{.5cm}

The cases covered by Theorem 4 do not represent the full unfolding
of the (Visible) Fold$-$Saddle singularity. Because of this, the
next two theorems are necessary. Each one of them describes a
distinct generic codimension two
singularity.\\

 \noindent {\bf Theorem 5.} {\it Take $\tau=v$ and $0 <\mu < \varepsilon_0$ in
Equation \eqref{eq fold-sela com parametros}. The
$(\lambda,\beta)-$plane contains essentially $13$ distinct typical
configurations representing $7$ distinct topological behaviors on
its bifurcation diagram (see Figure \ref{fig diagrama bif teo 4 5 e
6}).}

\vspace{.5cm}

 \noindent {\bf Theorem 6.} {\it Take $\tau=v$ and $-\varepsilon_0 <\mu < 0$
in Equation \eqref{eq fold-sela com parametros}. The
$(\lambda,\beta)-$plane contains essentially $13$ distinct typical
configurations representing $7$ distinct topological behaviors on
its bifurcation diagram (see Figure \ref{fig diagrama bif teo 4 5 e
6}).}

\vspace{.5cm}

Finally, we are able to state the main results of the paper.

\vspace{.5cm}

 \noindent {\bf Theorem A.} {\it  Equation \eqref{eq
fold-sela com parametros} with $\tau=i$ generically unfolds the
(Invisible) Fold$-$Saddle singularity. Moreover, its bifurcation
diagram exhibits $55$ distinct typical configurations representing
$11$ distinct topological behavior (see Figure \ref{fig diagrama
bifurcacao}).}

\vspace{.5cm}

 \noindent {\bf Theorem B.} {\it  Equation \eqref{eq
fold-sela com parametros} with $\tau=v$ generically unfolds the
(Visible) Fold$-$Saddle singularity. Moreover, its bifurcation
diagram exhibits $39$ distinct typical configurations representing
$21$ distinct topological behavior (see Figure \ref{fig diagrama
bifurcacao teo B}).}

\vspace{.5cm}

The paper is organized as follows: in Section 2 we give the basic
theory about Non$-$Smooth Vector Fields on the Plane, in Section 3
we prove Theorem 1, in Section 4 we prove Theorem 2, in Section 5 we
prove Theorem 3, in Section 6 we prove Theorem A and present the
Bifurcation Diagram of $\overline{Z}^{i}_{\lambda,\mu,\beta}$, in
Section 7 we prove Theorem 4, in Section 8 we prove Theorem 5, in
Section 9 we prove Theorem 6 and in Section 10 we prove Theorem B
and present the Bifurcation Diagram of
$\overline{Z}^{v}_{\lambda,\mu,\beta}$.


\section{Preliminaries}

We  distinguish the following regions on the discontinuity set
$\Sigma:$
\begin{itemize}
\item [(i)]$\Sigma_1\subseteq\Sigma$ is the \textit{sewing region} if
$(X.f)(Y.f)>0$ on $\Sigma_1$ .
\item [(ii)]$\Sigma_2\subseteq\Sigma$ is the \textit{escaping region} if
$(X.f)>0$ and $(Y.f)<0$ on $\Sigma_2$.
\item [(iii)]$\Sigma_3\subseteq\Sigma$ is the \textit{sliding region} if
$(X.f)<0$ and $(Y.f)>0$ on $\Sigma_3$.
\end{itemize}

\begin{figure}[!h]
\begin{center}
\psfrag{A}{$q$} \psfrag{B}{$q + Y(q)$} \psfrag{C}{$q + X(q)$}
\psfrag{D}{} \psfrag{E}{\hspace{1cm}$Z^\Sigma(q)$}
\psfrag{F}{\vspace{1cm}$\Sigma_2$} \psfrag{G}{} \epsfxsize=5.5cm
\epsfbox{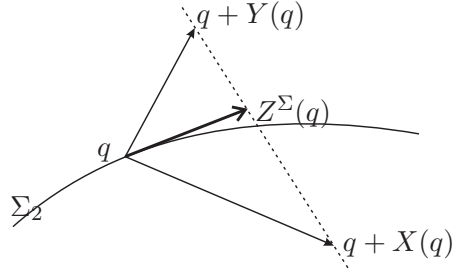} \caption{\small{Fillipov's
convention.}} \label{fig def filipov}
\end{center}
\end{figure}

Consider $Z \in \Omega^r.$ The \textit{sliding vector field}
associated to $Z$ is the vector field  $Z^s$ tangent to $\Sigma_3$
and defined at $q\in \Sigma_3$ by $Z^s(q)=m-q$ with $m$ being the
point where the segment joining $q+X(q)$ and $q+Y(q)$ is tangent to
$\Sigma_3$ (see Figure \ref{fig def filipov}). It is clear that if
$q\in \Sigma_3$ then $q\in \Sigma_2$ for $-Z$ and then we  can
define the {\it escaping vector field} on $\Sigma_2$ associated to
$Z$ by $Z^e=-(-Z)^s$. In what follows we use the notation $Z^\Sigma$
for both cases.\\

We say that $q\in\Sigma$ is a \textit{$\Sigma-$regular point} if
\begin{itemize}
\item [(i)] $(X.f(q))(Y.f(q))>0$
or
\item [(ii)] $(X.f(q))(Y.f(q))<0$ and $Z^{\Sigma}(q)\neq0$ (that is $q\in\Sigma_2\cup\Sigma_3$ and it is not a singular
point of $Z^{\Sigma}$).\end{itemize}

The points of $\Sigma$ which are not $\Sigma-$regular are called
\textit{$\Sigma-$singular}. We distinguish two subsets in the set of
$\Sigma-$singular points: $\Sigma^t$ and $\Sigma^p$. Any $q \in
\Sigma^p$ is called a \textit{pseudo equilibrium of $Z$} and it is
characterized by $Z^{\Sigma}(q)=0$. Any $q \in \Sigma^t$ is called a
\textit{tangential singularity} and is characterized by
$Z^{\Sigma}(q) \neq 0$ and
$X.f(q)Y.f(q) =0$ ($q$ is a contact point of $Z^{\Sigma}$).\\

%

A pseudo equilibrium $q \in \Sigma^p$ is a \textit{$\Sigma-$saddle}
provided one of the following condition is satisfied: (i)
$q\in\Sigma_2$ and $q$ is an attractor for $Z^{\Sigma}$ or (ii)
$q\in\Sigma_3$ and $q$ is a repeller for $Z^{\Sigma}$. A pseudo
equilibrium $q\in\Sigma^p$ is a $\Sigma-$\textit{repeller} (resp.
$\Sigma-$\textit{attractor}) provided $q\in\Sigma_2$ (resp. $q \in
\Sigma_3$) and $q$ is a repeller (resp. attractor) equilibrium point
for
$Z^{\Sigma}$. 

\begin{figure}[h!]
\begin{minipage}[b]{0.311\linewidth}
\begin{center}
\epsfxsize=3.3cm \epsfbox{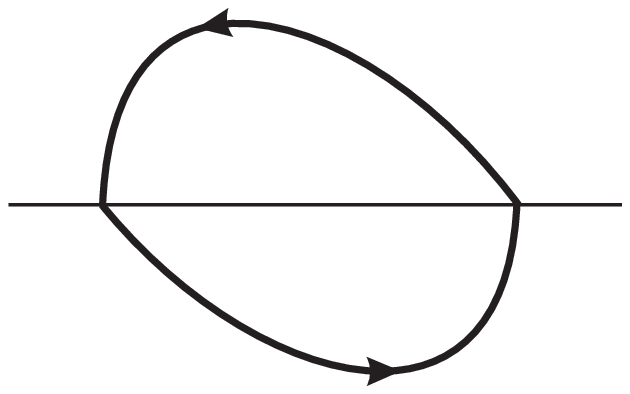} \caption{\small{Canard
cycle of kind I.}} \label{fig canard I}
\end{center}
\end{minipage} \hfill
\begin{minipage}[b]{0.311\linewidth}
\begin{center}
\psfrag{A}{$\Sigma=\Gamma$}\epsfxsize=2.6cm
\epsfbox{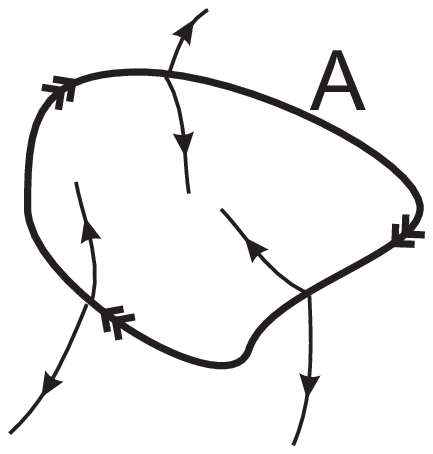} \caption{\small{Canard cycle of kind
II.}}\label{fig canard II}
\end{center}
\end{minipage} \hfill
\begin{minipage}[b]{0.36\linewidth}
\begin{center}
\epsfxsize=4cm  \epsfbox{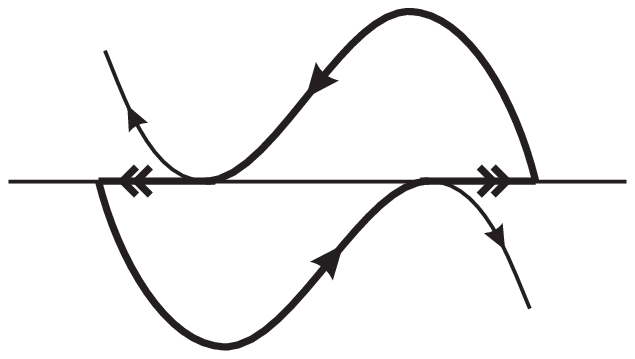}  \caption{\small{Canard
cycle of kind III.}}\label{fig canard}
\end{center}
\end{minipage}
\end{figure}

\begin{definition} Consider $Z \in \Omega^r.$
\begin{enumerate}
\item A curve $\Gamma$ is a \textbf{canard cycle} if
$\Gamma$ is closed and

  \begin{itemize}
  \item $\Gamma$ contains arcs of at least two of the vector fields $X |_{\Sigma_{+}}$, $Y |_{\Sigma_{-}}$ and $Z^{\Sigma}$ or is composed by a single arc of $Z^{\Sigma}$;

  \item the transition between arcs of $X$
  and arcs of $Y$ happens in sewing points;

  \item the transition between arcs of $X$
  (or $Y$) and arcs of $Z^{\Sigma}$ happens through
  $\Sigma-$fold points or regular points in the escape or sliding arc, respecting the orientation. Moreover
if $\Gamma\neq\Sigma$ then there exists at least one visible
$\Sigma-$fold point on each connected component of
$\Gamma\cap\Sigma$.
  \end{itemize}

\item Let $\Gamma$ be a canard cycle of $Z$. We say that

  \begin{itemize}
  \item $\Gamma$ is a \textbf{canard cycle of kind
  I} if $\Gamma$ meets $\Sigma$ just in sewing points;

  \item $\Gamma$ is a \textbf{canard cycle of kind
  II} se $\Gamma = \Sigma$;

  \item $\Gamma$ is a \textbf{canard cycle of kind
  III} if $\Gamma$ contains at least one visible $\Sigma-$fold point of $Z$.
  \end{itemize}

In Figures \ref{fig canard I}, \ref{fig canard II} and \ref{fig
canard} arise canard cycles of kind I, II and III respectively.

\item Let $\Gamma$ be  a canard cycle. We say that
$\Gamma$ is \textbf{hyperbolic} if

  \begin{itemize}
  \item $\Gamma$ is of kind I and $\eta'(p) \neq 1$,
  where $\eta$ is the first return map defined on a segment $T$ with $p\in T\pitchfork\gamma$;

  \item $\Gamma$ is of kind II;

  \item $\Gamma$ is of kind III and or $\Gamma\cap\Sigma\subseteq\Sigma_1\cup\Sigma_2$ or $\Gamma\cap\Sigma\subseteq\Sigma_1\cup\Sigma_3$.
  \end{itemize}
\end{enumerate}
\end{definition}

\begin{remark}The expression ``canard" is used here because these orbits are
limit periodic sets of singular perturbation problems (see
\cite{Eu-canard-cycles}).\end{remark}

\begin{definition}\label{centro}
Consider $Z \in \Omega^r$. A point $q \in \Sigma$ is a
\textbf{$\Sigma-$center} if there is a neighborhood $U$ of $q$ such
that an  one parameter family of canard cycles encircles $q$ and
foliates $U$.\\ 
\end{definition}

\begin{definition}\label{grafico} Consider $Z \in \Omega^r$. A closed path $\Delta$ is
a \textbf{$\Sigma-$graph} if it is a union of equilibria, pseudo
equilibria, tangential singularities of $Z$ and arcs of $Z$ joining
these points in such a way that $\Delta \cap \Sigma \neq \emptyset$.
Like for canard cycles, we say that $\Delta$ is a
\textbf{$\Sigma-$graph of kind I} if $\Delta \cap \Sigma \subset
\Sigma_1$, $\Delta$ is a \textbf{$\Sigma-$graph of kind II} if
$\Delta \cap \Sigma = \Delta$ and $\Delta$ is a
\textbf{$\Sigma-$graph of kind III} if $\Delta
\cap \Sigma \subsetneqq \Sigma_2 \cup \Sigma_3$.\\
\end{definition}

In what follows, in order to simplify the calculations, we take $\mu
=\alpha + 1$ in \eqref{eq fold-sela com parametros} and obtain the
following expression
\begin{equation}\label{eq fold-sela com parametros novos}
Z^{\tau}_{\lambda , \alpha, \beta} = \left\{
      \begin{array}{ll}
        X_{\lambda} = \left(
              \begin{array}{c}
                  1 \\
                \alpha_1(\tau)(x - \lambda)
\end{array}
      \right)
 & \hbox{if $y \geq 0$,} \\
         Y_{\alpha , \beta} = \left(
              \begin{array}{c}
                \frac{(1 + \alpha)}{2} x + \frac{(-1 + \alpha)}{2}(y + \beta) \\
                \frac{(-1 + \alpha)}{2} x + \frac{(1 + \alpha)}{2}(y + \beta)
\end{array}
      \right)& \hbox{if $y \leq 0$,}
      \end{array}
    \right.
\end{equation}
where   $\lambda$, $\beta \in (-1,1)$, $\alpha \in (-1 -
\varepsilon_0 ,-1 + \varepsilon_0)$, $\tau=i$ or $v$,
$\alpha_1(i)=-1$ and $\alpha_1(v)=1$. When it does not produce
confusion, in order to simplify the notation we
 use  $Z=(X,Y)$ or $Z_{\lambda , \alpha,
\beta} = (X , Y)$ instead $Z^{\tau}_{\lambda , \alpha, \beta} =
(X_{\lambda}, Y_{\alpha , \beta})$.\\

Given $Z=(X,Y)$, we describe some properties of both $X =
X_{\lambda}$ and $Y = Y_{\alpha, \beta}$.

The real number $\lambda$ measures how the $\Sigma-$fold point
$d=(\lambda,0)$ of $X$ is translated away from the origin. More
specifically, if $\lambda < 0$ then $d$ is translated to the left
hand side and if $\lambda
> 0$ then $d$ is
translated to the right hand side. 

\begin{figure}[!h]
\begin{minipage}[b]{0.485\linewidth}
\begin{center}\psfrag{A}{$\Sigma$} \psfrag{B}{$e$} \psfrag{G}{$Y.f=0$}
\psfrag{D}{$\Sigma$} \psfrag{E}{$X.f=0$} \psfrag{F}{$\gamma_1$}
\epsfxsize=4cm \epsfbox{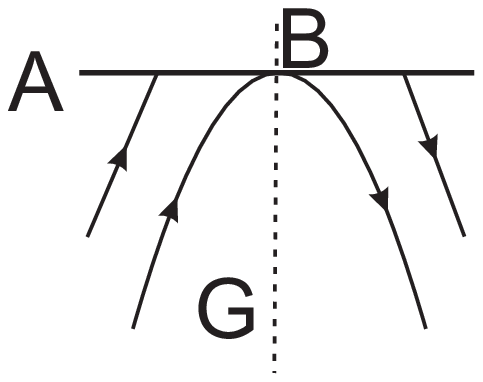} \caption{\small{Case
$Y^-$.}} \label{fig sela para cima}\end{center}
\end{minipage} \hfill
\begin{minipage}[b]{0.485\linewidth}
\begin{center}\psfrag{A}{$h$} \psfrag{B}{$i$} \psfrag{C}{$j$}
\psfrag{D}{$\sigma_2$} \psfrag{E}{$\sigma_1$} \psfrag{F}{$S$}
\psfrag{G}{$Y.f=0$} \epsfxsize=3cm \epsfbox{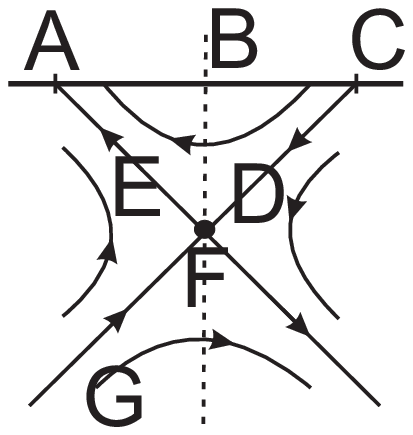}
\caption{\small{Case $Y^+$.}} \label{fig sela para baixo}
\end{center}\end{minipage} \hfill
\end{figure}

Some calculations show that the curve $Y.f=0$ is given by
$y=\frac{(1 - \alpha)}{(1 + \alpha)}x - \beta$. So the points of
this curve are equidistant from the separatrices when $\alpha = -1$.
It become closer to the stable separatrix of the saddle point $S =
S_{\alpha, \beta}$ when $\alpha \in (-1,-1+\varepsilon_0)$. It
become closer to the unstable separatrix of $S$ when $\alpha \in (-1
-\varepsilon_0,-1)$. Moreover, the smooth vector field $Y$ has
distinct types of contact with $\Sigma$ according with the
particular deformation considered. In this way, we have to consider
the following behaviors:

\begin{itemize}
\item $\mathbf{Y^- :}$ In this case $\beta < 0$. So  $S$ is translated to the
$y-$direction with $y>0$ (and $S$ is not visible for $Z$). It has a
visible $\Sigma-$fold point $e= e_{\alpha, \beta} = \Big( \frac{(1+
\alpha)}{(1 - \alpha)}\beta,0 \Big)$ (see Figure \ref{fig sela para
cima}).

\item $\mathbf{Y^0 :}$ In this case $\beta = 0$. So $S$ is not translated (see Figure \ref{fig fold-sela}).

\item $\mathbf{Y^+ :}$ In this case $\beta > 0$. So $S$ is translated to the
$y-$direction with $y<0$. It has an invisible $\Sigma-$fold point $i
= (i_1,i_2) = i_{\alpha, \beta} = \Big( \frac{(1+ \alpha)}{(1 -
\alpha)}\beta,0 \Big)$. Moreover, we distinguish two  points:  $h =
h_{\beta}=(-\beta,0)$ which is the intersection between the unstable
separatrix with $\Sigma$ and  $j = j_{\beta}=(\beta,0)$ which is the
intersection between the stable separatrix with $\Sigma$ (see Figure
\ref{fig sela para baixo}).
\end{itemize}

 In Figure \ref{fig sela para baixo} we
distinguish the arcs $\sigma_{1}$ of $Y$ joining the saddle point
$S$ of $Y$ to $h$ and  $\sigma_{2}$ of $Y$ joining $j$ to the saddle
point $S$ of $Y$.

\section{Proof of Theorem 1}\label{secao prova teorema 1}

\begin{figure}[!h]
\begin{center}
\psfrag{A}{$a$} \psfrag{B}{$b$} \psfrag{C}{$c$} \psfrag{D}{$\Sigma$}
\psfrag{E}{$X$} \psfrag{F}{$Y$} \psfrag{G}{$c+Y(c)$} \psfrag{H}{$c+
X(c)$} \psfrag{I}{$p(c)$} \epsfxsize=5cm
\epsfbox{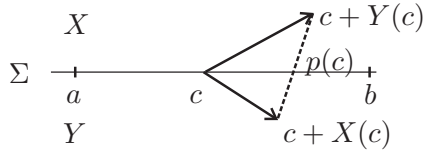} \caption{\small{Direction function.}}
\label{fig funcao direcao}
\end{center}
\end{figure}

In $(a,b)\subset \Sigma_2 \cup \Sigma_3$, consider  the point
$c=(c_{1},c_2)$, the vectors $X(c)=(d_1,d_2)$ and $Y(c)=(e_1,e_2)$
(as illustrated in Figure \ref{fig funcao direcao}). The straight
segment passing through $c+X(c)$ and $c + Y(c)$ meets $\Sigma$ in a
point $p(c)$. We define the C$^r-$map
$$
\begin{array}{cccc}
  p: & (a,b) & \longrightarrow & \Sigma \\
     & z & \longmapsto & p(z).
\end{array}
$$
We can choose local coordinates such that $\Sigma$ is the $x-$axis;
so $c=(c_1,0)$ and $p(c) \in \R \times \{ 0 \}$ can be identified
with points in $\R$. According with this identification, the
\textit{direction function} on $\Sigma$ is defined by
$$
\begin{array}{cccc}
  H: & (a,b) & \longrightarrow & \R \\
     & z & \longmapsto & p(z) - z.
\end{array}
$$

We obtain that $H$ is a C$^r-$map  and

\begin{itemize}
\item if $H(c) < 0$ then the orientation of $Z^{\Sigma}$ in a small neighborhood of $c$ is from $b$ to $a$;

\item if $H(c) = 0$ then  $c \in \Sigma^p$;

\item if $H(c) > 0$ then the orientation of $Z^{\Sigma}$ in a small neighborhood of $c$ is from $a$ to $b$.
\end{itemize}

Simple calculations show that $p(c_1) = \frac{e_2 (d_1+c_1)  -
d_2 (e_1+c_1)}{e_2 - d_2}$ and consequently,
\begin{equation}\label{eq H}
H(c_1) = \frac{e_2 d_1  - d_2 e_1}{e_2 - d_2}.
\end{equation}

We now in position to prove Theorem 1.

\begin{proof}[Proof of Theorem 1] In Cases $1_1$, $2_1$ and $3_1$ we assume
that $Y$ presents the behavior $Y^-$. In Cases $4_1$, $5_1$ and
$6_1$ we assume that $Y$ presents the behavior $Y^0$. In these cases
canard cycles are not allowed.

$\diamond$ \textit{Case $1_1$. $d<e$, Case $2_1$. $d=e$ and Case
$3_1$. $d>e$:} The points of $\Sigma$ outside the interval $(d,e)$
belong to $\Sigma_1$. The points inside this interval, when it is
not degenerated, belong to $\Sigma_3$ in Case $1_1$ and to
$\Sigma_{2}$ in Case $3_1$. In both cases $H(z)>0$ for all $z \in
(d,e)$. See Figure \ref{fig 1 teo 1}.

\begin{figure}[!h]
\begin{center}\psfrag{A}{\hspace{-.5cm}$\lambda < 0$} \psfrag{B}{$\lambda=0$} \psfrag{C}{$\lambda > 0$}\psfrag{0}{$1_1$} \psfrag{1}{$2_1$} \psfrag{2}{$3_1$}
\epsfxsize=8cm \epsfbox{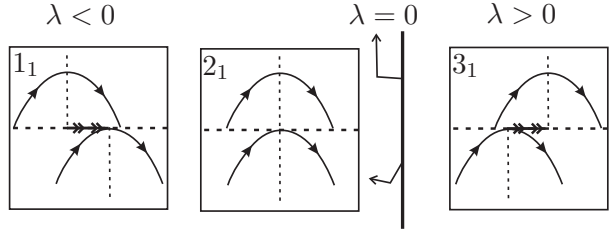} \caption{\small{Cases $1_1$,
$2_1$ and $3_1$.}} \label{fig 1 teo 1}
\end{center}
\end{figure}

$\diamond$ \textit{Case $4_1$. $d<S$, Case $5_1$. $d=S$ and Case
$6_1$. $d>S$:} The points of $\Sigma$ outside the interval $(d,S)$
belong to $\Sigma_1$. The points inside this interval, when it is
not degenerated, belong to $\Sigma_3$ in Case $4_1$ and to
$\Sigma_{2}$ in Case $6_1$. In both cases $H(z)>0$ for all $z \in
(d,S)$. See Figure \ref{fig 2 teo 1}.

\begin{figure}[!h]
\begin{center}\psfrag{A}{$\lambda < 0$} \psfrag{B}{$\lambda=0$} \psfrag{C}{$\lambda > 0$}\psfrag{0}{$4_1$} \psfrag{1}{$5_1$} \psfrag{2}{$6_1$}
\epsfxsize=8cm \epsfbox{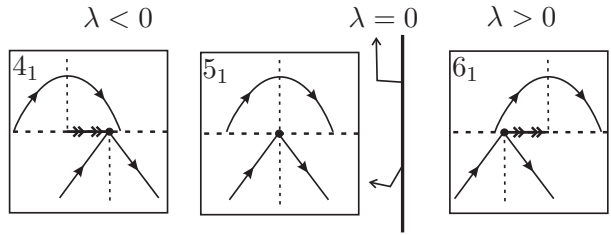} \caption{\small{Cases $4_1$,
$5_1$ and $6_1$.}} \label{fig 2 teo 1}
\end{center}
\end{figure}

In Cases $7_1 - 17_1$ we assume that $Y$ presents the behavior
$Y^+$.

$\diamond$ \textit{Case $7_1$. $\lambda< - \beta$, Case $8_1$.
$\lambda = - \beta$, Case $9_1$. $- \beta < \lambda < - \beta / 2$,
Case $10_1$. $\lambda = - \beta / 2$ and Case $11_1$. $- \beta / 2 <
\lambda < 0 $:} The points of $\Sigma$ outside the interval $(d,i)$
belong to $\Sigma_1$. The points inside this interval belong to
$\Sigma_3$. The direction function $H$ assumes positive values in a
neighborhood of $d$, negative values in a neighborhood of $i$ and
$H(\lambda \beta / (1 + \beta))=0$. So, by \eqref{eq H}, the
$\Sigma-$attractor $P=(\lambda \beta / (1 + \beta),0)$, nearby
$(0,0)$, is the unique pseudo equilibrium. In these cases canard
cycles are not allowed. See Figure \ref{fig 3a teo 1}.

\begin{figure}[!h]
\begin{center}\psfrag{A}{\hspace{-.8cm}$\lambda < - \beta$} \psfrag{B}{\hspace{-.4cm}$\lambda= - \beta$} \psfrag{C}{\hspace{-.35cm}$- \beta < \lambda
< -\beta / 2$} \psfrag{D}{$ \lambda = -\beta / 2$}
\psfrag{E}{\hspace{.35cm}$-\beta / 2 < \lambda <
0$}\psfrag{0}{$7_1$} \psfrag{1}{$8_1$} \psfrag{2}{$9_1$}
\psfrag{3}{$10_1$}\psfrag{4}{$11_1$} \epsfxsize=9cm
\epsfbox{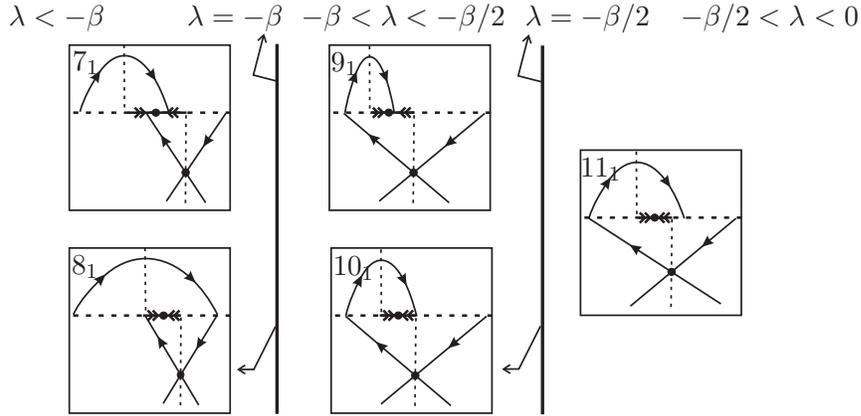} \caption{\small{Cases $7_1 - 11_1$.}}
\label{fig 3a teo 1}
\end{center}
\end{figure}

$\diamond$ \textit{Case $12_1$. $\lambda = 0$:}  Since $\alpha = -1$
and $d=i$, it is straightforward to show that each point $Q \in
(h,i)$ belongs to a closed curve composed by an arc of $X$ and an
arc of $Y$. So $d=i$ is a $\Sigma-$center. See Figure \ref{fig 3b
teo 1}.

\begin{figure}[!h]
\begin{center}\psfrag{4}{$12_1$} \epsfxsize=2.2cm
\epsfbox{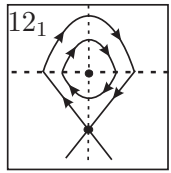} \caption{\small{Case $12_1$.}} \label{fig
3b teo 1}
\end{center}
\end{figure}

$\diamond$ \textit{Case $13_1$. $0<\lambda< \beta / 2$, Case $14_1$.
$\lambda = \beta / 2$, Case $15_1$. $\beta / 2 < \lambda < \beta$,
Case $16_1$. $\lambda = \beta$ and Case $17_1$. $\lambda
> \beta$:} The points of $\Sigma$ outside the interval $(i,d)$
belong to $\Sigma_1$ and the points inside this interval belong to
$\Sigma_2$. The direction function $H$ assumes positive values in a
neighborhood of $d$, negative values in a neighborhood of $i$ and
$H(\lambda \beta / (1 + \beta))=0$. So, by \eqref{eq H}, the
$\Sigma-$repeller $P=(\lambda \beta / (1 + \beta),0)$, nearby
$(0,0)$, is the unique pseudo equilibrium. In these cases canard
cycles are not allowed. See Figure \ref{fig 3c teo 1}.

\begin{figure}[!h]
\begin{center}\psfrag{A}{\hspace{-.7cm}$0<\lambda < \beta / 2$} \psfrag{B}{$\lambda= \beta / 2$} \psfrag{C}{\hspace{.2cm}$\beta / 2 < \lambda
< \beta$} \psfrag{D}{$ \lambda = \beta$} \psfrag{E}{$\lambda >
\beta$}\psfrag{0}{$13_1$} \psfrag{1}{$14_1$} \psfrag{2}{$15_1$}
\psfrag{3}{$16_1$}\psfrag{4}{$17_1$} \epsfxsize=9cm
\epsfbox{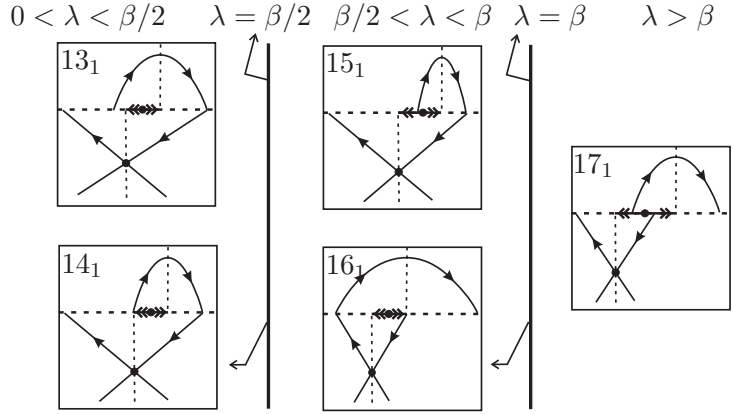} \caption{\small{Cases $13_1 - 17_1$.}}
\label{fig 3c teo 1}
\end{center}
\end{figure}

\begin{figure}[!h]
\begin{center}\psfrag{A}{$1_1$} \psfrag{B}{$2_1$} \psfrag{C}{$3_1$} \psfrag{D}{$4_1$} \psfrag{E}{$5_1$}
\psfrag{F}{$6_1$}  \psfrag{G}{$7_1$} \psfrag{H}{$8_1$}
\psfrag{I}{$9_1$} \psfrag{J}{$10_1$} \psfrag{K}{$11_1$}
\psfrag{L}{$12_1$}\psfrag{M}{$13_1$} \psfrag{N}{$14_1$}
\psfrag{O}{$15_1$} \psfrag{P}{$16_1$} \psfrag{Q}{$17_1$}
\epsfxsize=7cm \epsfbox{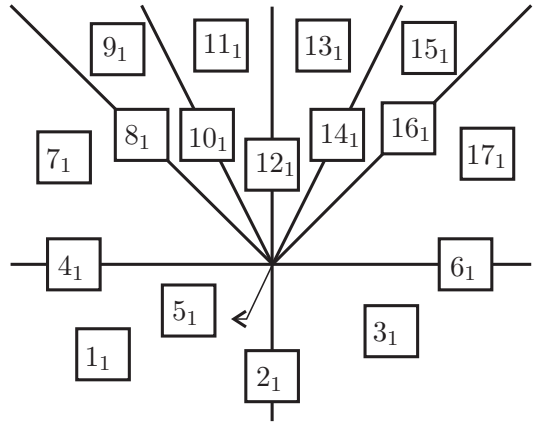}
\caption{\small{Bifurcation Diagram of Theorem 1.}} \label{fig
diagrama bif teo 1}
\end{center}
\end{figure}

The bifurcation diagram is illustrated in Figure \ref{fig diagrama
bif teo 1}.\end{proof}


\section{Proof of Theorem 2}\label{secao prova teorema 2}

\begin{proof}[Proof of Theorem 2] In Cases $1_2$, $2_2$ and $3_2$ we assume
that $Y$ presents the behavior $Y^-$.  In Cases $4_2$, $5_2$ and
$6_2$ we assume that $Y$ presents the behavior $Y^0$. In Cases $7_2
- 19_2$ we assume that $Y$ presents the behavior $Y^+$.

$\diamond$ \textit{Case $1_2$. $d<e$, Case $2_2$. $d=e$, Case $3_2$.
$d>e$, Case $4_2$. $d<S$, Case $5_2$. $d=S$ and Case $6_2$. $d>S$:}
Analogous to Cases $1_1$, $2_1$, $3_1$, $4_1$, $5_1$ and $6_1$.

$\diamond$ \textit{Case $7_2$. $\lambda< - \beta$, Case $8_2$.
$\lambda = - \beta$, Case $9_2$. $- \beta < \lambda < - \beta /
(1-\alpha)$, Case $10_2$. $\lambda = - \beta / (1-\alpha)$ and Case
$11_2$. $- \beta / (1-\alpha) < \lambda < 0 $:} Analogous to Cases
$7_1 - 11_1$ changing $- \beta / 2$ by $- \beta / (1-\alpha) = -
dist(h,i) / 2$, where $dist(h,i)$ is the distance between $h$ and
$i$. The unique pseudo equilibrium occurs in $P=(p^-,0)$ where
\begin{equation}\label{eq p1}\begin{array}{cl}p^-=& \dfrac{1}{2(\alpha+1)}
((1-\alpha)(1+\beta)+ \lambda(1+\alpha) +\\&-
\sqrt{((1-\alpha)(1+\beta)+ \lambda(1+\alpha))^2 - 4 \beta
(1+\alpha) (1+ \alpha + \lambda(1-\alpha))
}).\end{array}\end{equation}

$\diamond$ \textit{Case $12_2$. $\lambda = 0$:} The points of
$\Sigma$ outside the interval $(d,i)$ belong to $\Sigma_1$ and the
points inside this interval belong to $\Sigma_3$. The direction
function $H$ assumes positive values in a neighborhood of $d$,
negative values  in a  neighborhood of $i$ and $H(p^+_{0},0)=0$
where $p^+_0$ is given by \eqref{eq p1} with $\lambda=0$. So
$P=(p^+_0,0)$ is a $\Sigma-$attractor. Since
 $e=0$, it is easy to see that there is an arc
$\gamma_{1}^X$ of $X$ connecting the points $h$ and $j$. It
generates a $\Sigma-$graph $\Gamma = \gamma_{1}^{X} \cup \sigma_2
\cup S \cup \sigma_1$ of kind I. Since $-1 < \alpha < -1 +
\varepsilon_0$, it is straight forward to show that the
\textit{First Return Map} $\eta=\varphi_{Y}\circ \varphi_{X}$, where
$$
\begin{array}{cccc}
  \varphi_{X}: & \Sigma & \rightarrow & \Sigma \\
     & z=(x,0) & \longmapsto & (-x+2 \lambda, 0)
\end{array}\mbox{ and }$$$$
\begin{array}{cccc}
  \varphi_{Y}: & (i,j)\subset \Sigma & \rightarrow & (h,i)\subset \Sigma  \\
     & z=(x,0) & \longmapsto & \Big( \dfrac{x(i_1+\beta)-2i_{1}^{2}}{\beta -i_1}, 0
     \Big)
\end{array},$$has derivative bigger than $1$ in the interval $(h,d)$. By
consequence, $\Gamma$ is a repeller for the trajectories inside it
and in this case canard cycles are not allowed. See Figure \ref{fig
1 teo 2}.

\begin{figure}[!h]
\begin{center}\psfrag{A}{$\lambda =0$} \psfrag{B}{$0<\lambda < i_1$} \psfrag{C}{$\lambda = i_1$} \psfrag{D}{$ \lambda = -\frac{\beta}{2}$}
\psfrag{E}{$-\frac{\beta}{2} < \lambda < 0$}\psfrag{0}{$12_2$}
\psfrag{1}{$13_2$} \psfrag{2}{$14_2$}
\psfrag{3}{$10_1$}\psfrag{4}{$11_1$} \epsfxsize=10cm
\epsfbox{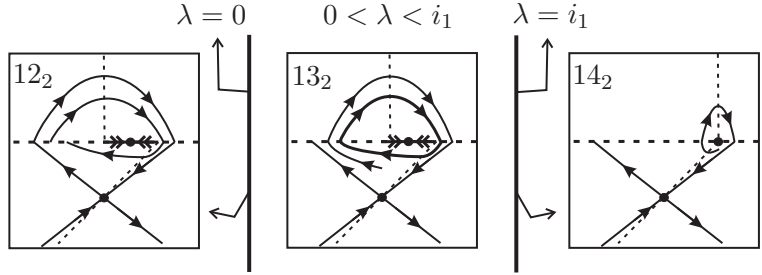} \caption{\small{Cases $12_2$, $13_2$ and
$13_2$.}} \label{fig 1 teo 2}
\end{center}
\end{figure}

$\diamond$ \textit{Case $13_2$. $0<\lambda < i_1$:} The distribution
of the connected components of $\Sigma$ and the behavior of $H$ are
the same of Case $12_2$ with $P=(p^+_\lambda,0)$ where $p^+_\lambda$
is given by \eqref{eq p1}. Since $0<\lambda < i_1$, there is an arc
$\gamma_{1}^X$ of $X$ connecting the point $j$ to a point $k_1 \in
(h,d)$. Also there is an arc $\gamma_{1}^{Y}$ of $Y$ connecting the
point $k_1$ to a point $l_1 \in (i,j)$. Repeating this argument, we
can find an increasing sequence $(k_i)_{i \in \N}$. We can prove
that there is an interval $I\subset(k_1,d)$ such that
$\eta'=(\varphi_{Y}\circ \varphi_{X})' < 1$. As $P$ is a
$\Sigma-$attractor, there is an interval $J \subset (k_1,d)$ such
that $\eta' > 1$. Moreover, there exists an unique point $Q \in
(k_1,d)$ given by $Q=((-i_{1}^{2}+\lambda (i_1 + \beta))/\beta,0)$
such that $\eta'=1$. By $Q$ passes a repeller canard cycle $\Gamma$
of kind I. See Figure \ref{fig 1 teo 2}.

$\diamond$ \textit{Case $14_2$. $\lambda = i_1$:} Every point of
$\Sigma$ belongs to $\Sigma_1$ except the point $d=i$. As in the
previous case, we can construct sequences $(k_i)_{i\in \N}$ and
$(l_i)_{i\in \N}$. Since $e=i_1$, we have that $k_i \rightarrow d$
and $l_i \rightarrow d$. So $d$ is a non generic tangential
singularity of repeller kind. In this case canard cycles are not
allowed. See Figure \ref{fig 1 teo 2}.

$\diamond$ \textit{Case $15_2$. $i_1<\lambda< \alpha \beta /
(1-\alpha)$, Case $16_2$. $\lambda = \alpha \beta / (1-\alpha)$,
Case $17_2$. $\alpha \beta / (1-\alpha) < \lambda < \beta$, Case
$18_2$. $\lambda = \beta$ and Case $19_2$. $\lambda
> \beta$:}  Analogous to
Cases $13_1 - 17_1$ changing $\beta / 2$ by $ \alpha \beta /
(1-\alpha)= - dist(i,j) / 2$. The unique pseudo equilibrium occurs
in $P=(p^-,0)$ where $p^-$ is given by \eqref{eq p1}.

\begin{figure}[!h]
\begin{center}\psfrag{A}{$1_2$} \psfrag{B}{$2_2$} \psfrag{C}{$3_2$} \psfrag{D}{$4_2$} \psfrag{E}{$5_2$}
\psfrag{F}{$6_2$}  \psfrag{G}{$7_2$} \psfrag{H}{$8_2$}
\psfrag{I}{$9_2$} \psfrag{J}{$10_2$} \psfrag{K}{$11_2$}
\psfrag{L}{$12_2$}\psfrag{M}{$13_2$} \psfrag{N}{$14_2$}
\psfrag{O}{$15_2$} \psfrag{P}{$16_2$}
\psfrag{Q}{$17_2$}\psfrag{R}{$18_2$} \psfrag{S}{$19_2$}
\epsfxsize=10cm \epsfbox{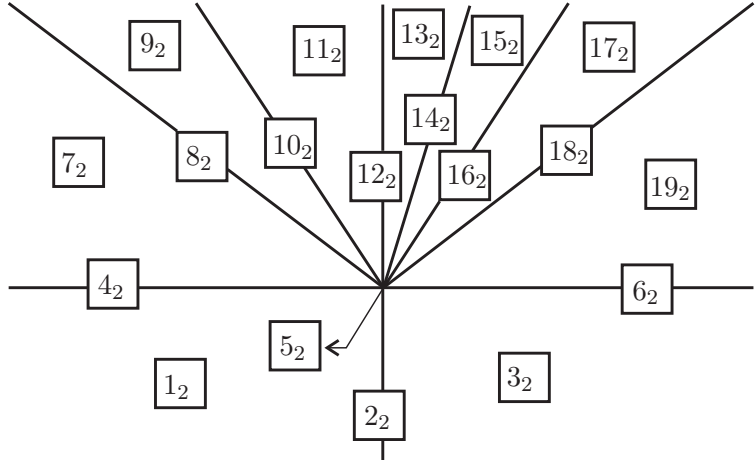}
\caption{\small{Bifurcation Diagram of Theorem 2.}} \label{fig
diagrama bif teo 2}
\end{center}
\end{figure}

The bifurcation diagram is illustrated in Figure \ref{fig diagrama
bif teo 2}.\end{proof}


\section{Proof of Theorem 3}\label{secao prova teorema 3}

\begin{proof}[Proof of Theorem 3] In Cases $1_3$, $2_3$ and $3_3$ we assume
that $Y$ presents the behavior $Y^-$.  In Cases $4_3$, $5_3$ and
$6_3$ we assume that $Y$ presents the behavior $Y^0$. In Cases $7_3
- 19_3$ we assume that $Y$ presents the behavior $Y^+$.

$\diamond$ \textit{Case $1_3$. $d<e$, Case $2_3$. $d=e$, Case $3_3$.
$d>e$, Case $4_3$. $d<S$, Case $5_3$. $d=S$ and Case $6_3$. $d>S$:}
Analogous to Cases $1_1$, $2_1$, $3_1$, $4_1$, $5_1$ and $6_1$.

$\diamond$ \textit{Case $7_3$. $\lambda< - \beta$, Case $8_3$.
$\lambda = - \beta$, Case $9_3$. $- \beta < \lambda < - \beta /
(1-\alpha)$, Case $10_3$. $\lambda = - \beta / (1-\alpha)$ and Case
$11_3$. $- \beta / (1-\alpha) < \lambda < i_1 $:} Analogous to Cases
$7_1 - 11_1$ changing $- \beta / 2$ by $- \beta / (1-\alpha) = -
dist(h,i) / 2$.  The unique pseudo equilibrium occurs in
$P=(p^{-},0)$ where $p^-$ is given by \eqref{eq p1}.

$\diamond$ \textit{Case $12_3$. $\lambda = i_1$:} Analogous to Case
$14_2$ except that here $d$ is an attractor, i.e., there is a change
of stability. See Figure \ref{fig 1 teo 3}.

\begin{figure}[!h]
\begin{center}\psfrag{A}{$\lambda = i_1$} \psfrag{B}{$i_1<\lambda < 0$} \psfrag{C}{$\lambda =0$}
\psfrag{1}{$13_3$} \psfrag{2}{$14_3$} \psfrag{0}{$12_3$}
\epsfxsize=10cm \epsfbox{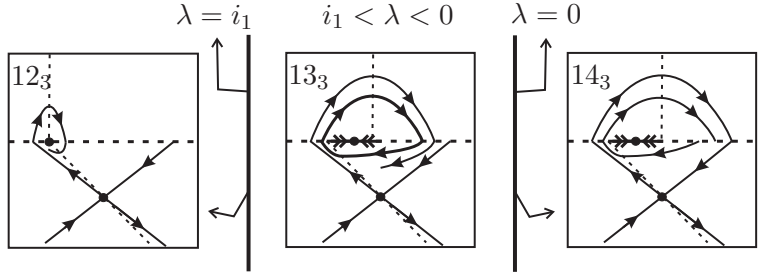} \caption{\small{Cases
$12_3$, $13_3$ and $13_3$.}} \label{fig 1 teo 3}
\end{center}
\end{figure}

$\diamond$ \textit{Case $13_3$. $i_1 <\lambda < 0$:} Analogous to
Case $13_2$ except that there is a change of stability on
$P=(p^{-},0)$
, which is a $\Sigma-$repeller, and on $\Gamma$, which is an
attractor canard cycle of kind I. See Figure \ref{fig 1 teo 3}.

$\diamond$ \textit{Case $14_3$. $\lambda = 0$:} Analogous to Case
$12_2$ except that occurs a change of stability on $P=(p^{-},0)$
, which is a
$\Sigma-$repeller, and on $\Gamma$, which is an attractor for the
trajectories inside it. See Figure \ref{fig 1 teo 3}.

$\diamond$ \textit{Case $15_3$. $0<\lambda< \alpha \beta /
(1-\alpha)$, Case $16_3$. $\lambda = \alpha \beta / (1-\alpha)$,
Case $17_3$. $\alpha \beta / (1-\alpha) < \lambda < \beta$, Case
$18_3$. $\lambda = \beta$ and Case $19_2$. $\lambda
> \beta$:}  Analogous to
Cases $13_1 - 17_1$ changing $\beta / 2$ by $\alpha \beta /
(1-\alpha)= - dist(i,j) / 2$. The unique pseudo equilibrium occurs
in $P=(p^{-},0)$.

\begin{figure}[!h]
\begin{center}\psfrag{A}{$1_3$} \psfrag{B}{$2_3$} \psfrag{C}{$3_3$} \psfrag{D}{$4_3$} \psfrag{E}{$5_3$}
\psfrag{F}{$6_3$}  \psfrag{G}{$7_3$} \psfrag{H}{$8_3$}
\psfrag{I}{$9_3$} \psfrag{J}{$10_3$} \psfrag{K}{$11_3$}
\psfrag{L}{$12_3$}\psfrag{M}{$13_3$} \psfrag{N}{$14_3$}
\psfrag{O}{$15_3$} \psfrag{P}{$16_3$}
\psfrag{Q}{$17_3$}\psfrag{R}{$18_3$} \psfrag{S}{$19_3$}
\epsfxsize=10cm \epsfbox{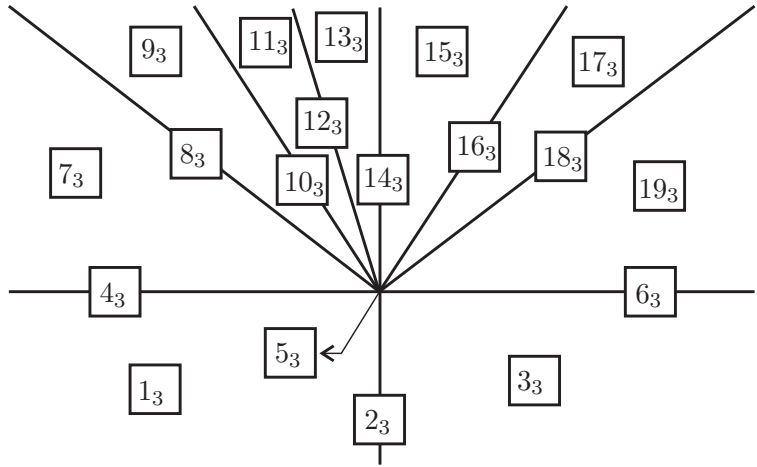}
\caption{\small{Bifurcation Diagram of Theorem 3.}} \label{fig
diagrama bif teo 3}
\end{center}
\end{figure}

The bifurcation diagram is illustrated in Figure \ref{fig diagrama
bif teo 3}.\end{proof}


\section{Proof of Theorem A}\label{secao prova teorema A}

\begin{proof}[Proof of Theorem A] Since in Equation \eqref{eq fold-sela com
parametros novos} we can take $\alpha$ in the interval
$(-\infty,0)$, from Theorems 1, 2 and 3 we derive that this
equation, with $\tau=i$, unfolds generically the (Invisible)
Fold$-$Saddle singularity.

Observe that the bifurcation diagram contain all the typical
configurations and all the distinct topological behavior described
in Theorems 1, 2 and 3. So, the number of typical configurations is
$55$ and the number of distinct topological behaviors is $11$.
Moreover, each topological behavior can be represented respectively
by the Cases $1_1$, $4_1$,  $7_1$, $12_1$, $13_1$, $12_2$, $13_2$,
$14_2$, $12_3$, $13_3$ and $14_3$.

The full behavior of the three$-$parameter family of non$-$smooth
vector fields presenting the normal form \eqref{eq fold-sela com
parametros novos}, with $\tau=i$, is illustrated in Figure \ref{fig
diagrama bifurcacao}  where we consider a sphere around the point
$(\lambda , \mu, \beta) = (0,0,0)$ with a small ray and so we make a
stereographic projection defined on the entire sphere, except the
south pole. Still in relation with this figure, the numbers pictured
correspond to the occurrence of the cases described in the previous
theorems. As expected, the cases $5_1$ and $5_2$ are not represented
in this figure because they are, respectively, the center and the
south pole of the sphere.\end{proof}

\begin{figure}[h!]
\begin{center}
\psfrag{A}{$1_2$}\psfrag{B}{$1_1$}\psfrag{C}{$1_3$}
\psfrag{D}{$2_3$}\psfrag{E}{$2_2$}
\psfrag{F}{$2_1$}\psfrag{G}{$3_2$}\psfrag{H}{$3_1$}
\psfrag{I}{$3_3$}\psfrag{J}{$4_3$}
\psfrag{K}{$5_3$}\psfrag{L}{$4_1$}\psfrag{M}{$4_2$}\psfrag{N}{$6_1$}
\psfrag{O}{$6_3$}
\psfrag{P}{$6_2$}\psfrag{Q}{$7_2$}\psfrag{R}{$7_1$}\psfrag{S}{$7_3$}
\psfrag{T}{$8_2$}
\psfrag{U}{$8_1$}\psfrag{V}{$8_3$}\psfrag{X}{$9_2$}\psfrag{Y}{$10_1$}
\psfrag{W}{$10_2$}
\psfrag{Z}{$9_1$}\psfrag{0}{$11_2$}\psfrag{1}{$10_3$}\psfrag{2}{$9_3$}
\psfrag{3}{$12_3$}
\psfrag{4}{$11_1$}\psfrag{5}{$11_3$}\psfrag{6}{$12_2$}\psfrag{7}{$12_1$}
\psfrag{8}{$13_3$} \psfrag{9}{$13_2$}\psfrag{`}{$\hspace{-.8cm}-1<
\alpha < 0$}\psfrag{~}{$\alpha =
-1$}\psfrag{:}{$17_2$}\psfrag{@}{$18_2$}
\psfrag{#}{$14_1$}\psfrag{^}{$17_3$}\psfrag{&}{$17_1$}\psfrag{*}{$18_3$}\psfrag{}{$$}
\psfrag{-}{$14_3$}\psfrag{_}{$16_1$}\psfrag{=}{$19_2$}\psfrag{+}{$19_3$}\psfrag{}{$$}
\psfrag{|}{$\lambda
>0$}\psfrag{;}{$15_1$}\psfrag{,}{$\hspace{.4cm}\alpha<-1$}\psfrag{<}{$15_2$}
\psfrag{.}{$15_3$} \psfrag{>}{$13_1$}\psfrag{?}{$\beta <
0$}\psfrag{a}{$\lambda =
0$}\psfrag{b}{$\beta=0$}\psfrag{c}{$\beta>0$} \psfrag{d}{$\lambda <
0$}\psfrag{e}{$16_2$}\psfrag{f}{$16_3$}\psfrag{g}{$14_2$}\psfrag{}{$$}
\epsfxsize=13cm \epsfbox{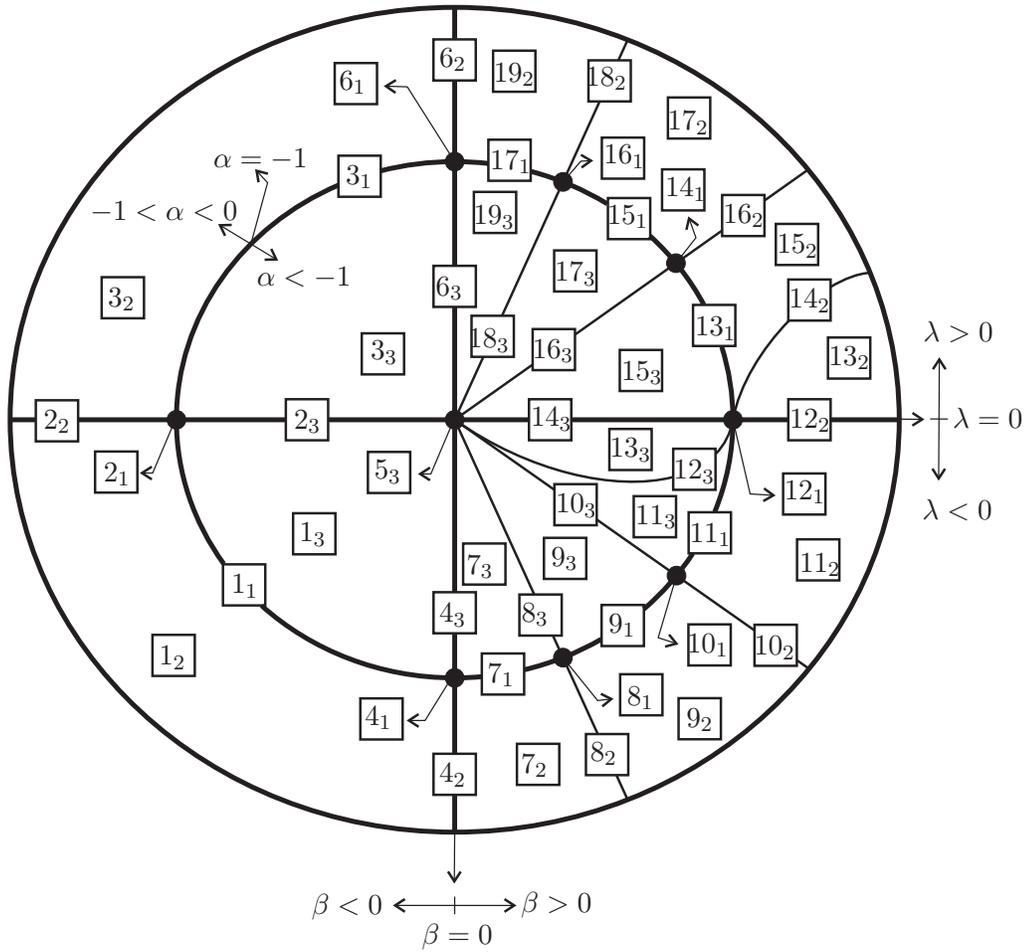}
\caption{\small{Bifurcation diagram of the (Invisible) Fold$-$Saddle
singularity.}} \label{fig diagrama bifurcacao}
\end{center}
\end{figure}

\section{Proof of Theorem 4}\label{secao prova teorema 4}


\begin{proof}[Proof of Theorem 4] Since $X$ has a unique $\Sigma-$fold point which is visible we conclude that canard cycles are
not allowed.

In Cases $1_4$, $2_4$ and $3_4$ we assume that $Y$ presents the
behavior $Y^-$. In Cases $4_4$, $5_4$ and $6_4$ we assume that $Y$
presents the behavior $Y^0$. In these cases, when
 it is well defined, the direction function $H$ assumes positive values.

$\diamond$ \textit{Case $1_4$. $d<e$:} The points of $\Sigma$ inside
the interval $(d,e)$ belong to $\Sigma_1$. The points on the left of
 $d$ belong to $\Sigma_3$ and the points on the right of
 $e$ belong to $\Sigma_2$. See Figure \ref{fig 1 teo 4}.

$\diamond$ \textit{Case $2_4$. $d=e$:} Here $\Sigma_1 = \emptyset$.
  The vector fields $X$ and $Y$ are linearly
 dependent on $d=e$ which is a tangential singularity. Moreover, it is an attractor
 for the trajectories of $Z$ crossing $\Sigma_3$ and  a repeller for the trajectories of $Z$ crossing
 $\Sigma_2$. See Figure \ref{fig 1 teo 4}.

$\diamond$ \textit{Case $3_4$. $d>e$:} The points of $\Sigma$ inside
the interval $(e,d)$ belong to $\Sigma_1$. The points on the left of
 $e$ belong to $\Sigma_3$ and the points on the right of
 $d$ belong to $\Sigma_2$. See Figure \ref{fig 1 teo 4}.

\begin{figure}[!h]
\begin{center}\psfrag{A}{$\lambda < 0$} \psfrag{B}{$\lambda=0$} \psfrag{C}{$\lambda > 0$}\psfrag{0}{$1_4$} \psfrag{1}{$2_4$} \psfrag{2}{$3_4$}
\epsfxsize=8cm \epsfbox{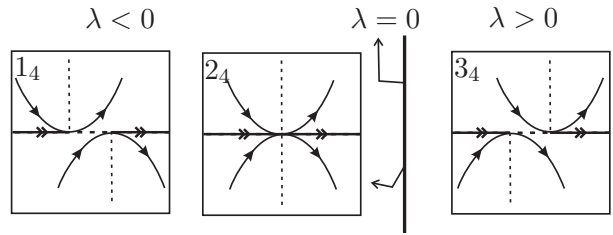} \caption{\small{Cases $1_4$,
$2_4$ and $3_4$.}} \label{fig 1 teo 4}
\end{center}
\end{figure}


$\diamond$ \textit{Case $4_4$. $d<S$:} The points of $\Sigma$ inside
the interval $(d,S)$ belong to $\Sigma_1$. The points on the left of
 $d$ belong to $\Sigma_3$ and the points on the right of
 $S$ belong to $\Sigma_2$. See Figure \ref{fig 2 teo 4}.

$\diamond$ \textit{Case $5_4$. $d=S$:} Here $\Sigma_1=\emptyset$ and
$S$ is an attractor
 for the trajectories of $Z$ crossing $\Sigma_3$ and it is a repeller for the trajectories of $Z$ crossing
 $\Sigma_2$. See Figure \ref{fig 2 teo 4}.

$\diamond$ \textit{Case $6_4$. $d>S$:} The points of $\Sigma$ inside
the interval $(d,S)$ belong to $\Sigma_1$. The points on the left of
 $S$ belong to $\Sigma_3$ and the points on the right of
 $d$ belong to $\Sigma_2$. See Figure \ref{fig 2 teo 4}.

\begin{figure}[!h]
\begin{center}\psfrag{A}{\hspace{-.5cm}$\lambda < 0$} \psfrag{B}{$\lambda=0$} \psfrag{C}{$\lambda > 0$}\psfrag{0}{$4_4$} \psfrag{1}{$5_4$} \psfrag{2}{$6_4$}
\epsfxsize=8cm \epsfbox{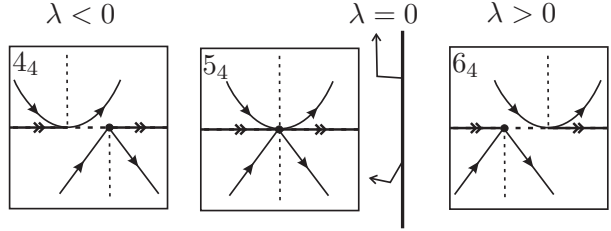} \caption{\small{Cases $4_4$,
$5_4$ and $6_4$.}} \label{fig 2 teo 4}
\end{center}
\end{figure}

In Cases $7_4 - 13_4$ we assume that $Y$ presents the behavior
$Y^+$.

$\diamond$ \textit{Case $7_4$. $d < h$, Case $8_4$. $d = h$ and Case
$9_4$. $h < d < i$:} The points of $\Sigma$ inside the interval
$(d,i)$ belong to $\Sigma_1$. The points on the left of
 $d$ belong to $\Sigma_3$ and the points on the right of
 $i$ belong to $\Sigma_2$. The direction function $H$ assumes positive
values on $\Sigma_3$ and negative values in a neighborhood of $i$.
Moreover, $H(\beta \lambda / (-1 + \beta))=0$ and the
$\Sigma-$repeller $P=(\beta \lambda / (-1 + \beta),0)$ is the unique
pseudo equilibrium. See Figure \ref{fig 3a teo 4}.

$\diamond$ \textit{Case $10_4$. $d=i$:} Here $\Sigma_1 = \emptyset$.
The vector fields $X$ and $Y$ are linearly
 dependent on the tangential singularity $d=i$.
A straightforward calculation shows that $H(z)=(1-\beta)/2 \neq 0$
for all $z \in \Sigma/\{ d \}$. So $d=i$ is an attractor
 for the trajectories of $Z$ crossing $\Sigma_3$ and  a repeller for the trajectories of $Z$ crossing
 $\Sigma_2$. Moreover, $\Delta=\{ d \} \cup \overline{d j} \cup \sigma_2 \cup \{ S \} \cup \sigma_1 \cup \overline{h d}$
 is a $\Sigma-$graph of kind III in such a way that each $Q$ in its interior belongs to another $\Sigma-$graph of kind III passing through $d$.  See Figure
\ref{fig 3a teo 4}.

\begin{figure}[!h]
\begin{center}\psfrag{A}{$\lambda < -\beta$} \psfrag{B}{$\lambda = - \beta$} \psfrag{C}{$- \beta < \lambda
< 0$} \psfrag{D}{$ \lambda = 0$}\psfrag{0}{$7_4$} \psfrag{1}{$8_4$}
\psfrag{2}{$9_4$} \psfrag{3}{$10_4$} \epsfxsize=11.5cm
\epsfbox{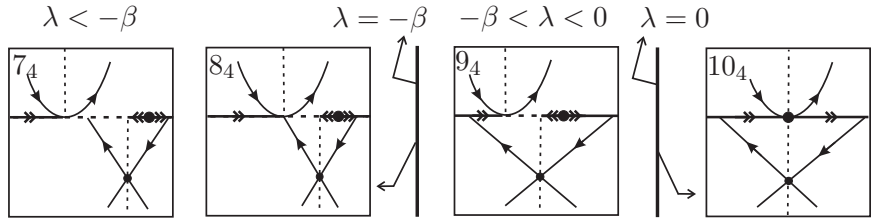} \caption{\small{Cases $7_4 - 10_4$.}}
\label{fig 3a teo 4}
\end{center}
\end{figure}

$\diamond$ \textit{Case $11_4$. $i < d < j$, Case $12_4$. $d = j$
and Case $13_4$. $j < d$:} The points of $\Sigma$ inside the
interval $(i,d)$ belong to $\Sigma_1$. The points on the left of
 $i$ belong to $\Sigma_3$ and the points on the right of
 $d$ belong to $\Sigma_2$. The direction function $H$ assumes
 positive
values on $\Sigma_2$ and negative values in a neighborhood of $i$.
Moreover, $H(\beta \lambda / (-1 + \beta))=0$ and the
$\Sigma-$attractor $P=(\beta \lambda / (-1 + \beta),0)$ is the
unique pseudo equilibrium. See Figure \ref{fig 3b teo 4}.

\begin{figure}[!h]
\begin{center}\psfrag{A}{$0<\lambda < \beta$} \psfrag{B}{$\lambda = \beta$} \psfrag{C}{$\beta < \lambda$}
\psfrag{0}{$11_4$} \psfrag{1}{$12_4$} \psfrag{2}{$13_4$}
 \epsfxsize=8cm
\epsfbox{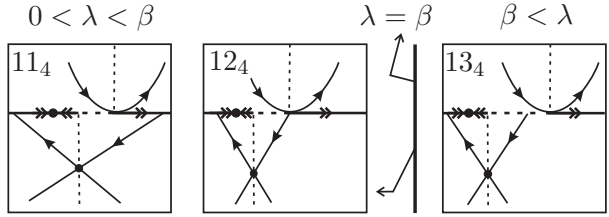} \caption{\small{Cases $11_4 - 13_4$.}}
\label{fig 3b teo 4}
\end{center}
\end{figure}

\begin{figure}[!h]
\begin{center}\psfrag{A}{\hspace{.15cm}$1$} \psfrag{B}{\hspace{.1cm}$2$}
 \psfrag{C}{\hspace{.15cm}$3$} \psfrag{D}{\hspace{.12cm}$4$} \psfrag{E}{\hspace{.13cm}$5$}
\psfrag{F}{\hspace{.1cm}$6$}  \psfrag{G}{\hspace{.15cm}$7$}
\psfrag{H}{\hspace{.15cm}$8$} \psfrag{I}{\hspace{.17cm}$9$}
\psfrag{J}{\hspace{.15cm}$10$} \psfrag{K}{\hspace{.15cm}$11$}
\psfrag{L}{\hspace{.1cm}$12$}\psfrag{M}{\hspace{.1cm}$13$}
\epsfxsize=7cm \epsfbox{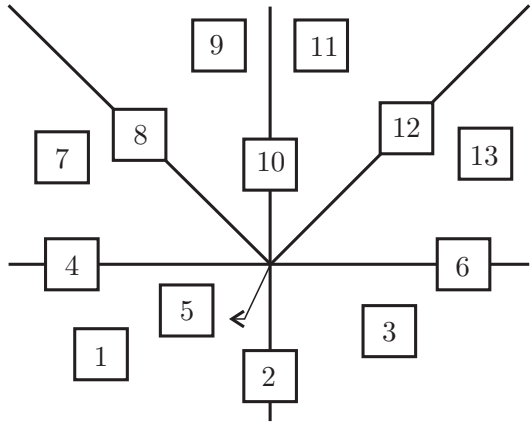}
\caption{\small{Bifurcation Diagram of Theorems 4, 5 and 6.}}
\label{fig diagrama bif teo 4 5 e 6}
\end{center}
\end{figure}

The bifurcation diagram is illustrated in Figure \ref{fig diagrama
bif teo 4 5 e 6}. Each topological behavior can be represented
respectively by Cases $1_4$, $2_4$, $4_4$, $5_4$, $7_4$, $10_4$ and
$11_4$.\end{proof}


\section{Proof of Theorem 5}\label{secao prova teorema 5}

\begin{proof}[Proof of Theorem 5] The
direction function $H$ has a  root $Q=(q,0)$
where\begin{equation}\label{eq p novo}\begin{array}{cl}q=&
\dfrac{1}{2(\alpha+1)} ((-1+\alpha)(1-\beta)- \lambda(1+\alpha) +\\&
+ \sqrt{((-1+\alpha)(1-\beta)- \lambda(1+\alpha))^2 + 4 \beta
(1+\alpha) (1+ \alpha + \lambda(-1+\alpha))
}).\end{array}\end{equation}Moreover, $H$ assumes positive values on
the right of $Q$ and negative values on the left of $Q$. Note that
when $\alpha \rightarrow -1$ so $Q \rightarrow - \infty$ under the
line $\{y=0 \}$ and it occurs the configurations showed in Theorem
4.

In Cases $1_5$, $2_5$ and $3_5$ we assume that $Y$ presents the
behavior $Y^-$.  In Cases $4_5$, $5_5$ and $6_5$ we assume that $Y$
presents the behavior $Y^0$. In Cases $7_5 - 13_5$ we assume that
$Y$ presents the behavior $Y^+$.

$\diamond$ \textit{Case $1_5$. $d<e$, Case $2_5$. $d=e$, Case $3_5$.
$d>e$, Case $4_5$. $d<S$, Case $5_5$. $d=S$ and Case $6_5$. $d>S$:}
Analogous to Cases $1_4$, $2_4$, $3_4$, $4_4$, $5_4$ and $6_4$
respectively, except that here it appears the $\Sigma-$saddle $Q$ on
the left of $d$ and $e$ or $S$. See Figure \ref{fig 1 teo 5}.

\begin{figure}[!h]
\begin{center}\psfrag{A}{\hspace{-2.6cm}$\lambda < (1+\alpha) \beta/(1-\alpha)$}
\psfrag{B}{\hspace{-2.2cm}$\lambda=(1+\alpha) \beta/(1-\alpha)$}
\psfrag{C}{$\lambda > (1+\alpha) \beta/(1-\alpha)$}
\psfrag{0}{$1_5$} \psfrag{1}{$2_5$} \psfrag{2}{$3_5$} \epsfxsize=8cm
\epsfbox{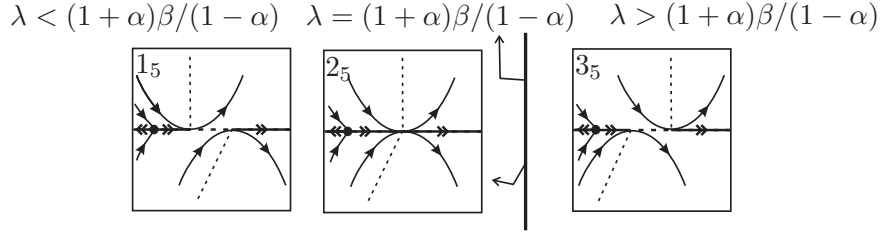} \caption{\small{Cases $1_5$, $2_5$ and
$3_5$.}} \label{fig 1 teo 5}
\end{center}
\end{figure}

%


$\diamond$ \textit{Case $7_5$. $d< h$, Case $8_5$. $d = h$, Case
$9_5$. $h < d < i$:} Analogous to Cases $7_4 - 9_4$, except that
here it appears the $\Sigma-$saddle $Q$ on the left of $d$ and $i$.
Here $P=(p,0)$ where\begin{equation}\label{eq p novo
2}\begin{array}{cl}p=& \dfrac{1}{2(\alpha+1)} ((-1+\alpha)(1-\beta)-
\lambda(1+\alpha) +\\& - \sqrt{((-1+\alpha)(1-\beta)-
\lambda(1+\alpha))^2 + 4 \beta (1+\alpha) (1+ \alpha +
\lambda(-1+\alpha)) }).\end{array}\end{equation}

$\diamond$ \textit{Case $10_5$. $d=i$:} Analogous to Case $10_4$,
except that here appear the $\Sigma-$saddle $Q$ on the left of
$d=i$.

$\diamond$ \textit{Case $11_5$. $i < d < j$, Case $12_5$. $d = j$
and Case $13_5$. $j < d$:}  Analogous to Cases $11_4 - 13_4$, except
that here it appears the $\Sigma-$saddle $Q$ on the left of $d$ and
$i$.

%

The bifurcation diagram is illustrated in Figure \ref{fig diagrama
bif teo 4 5 e 6}. Each topological behavior can be represented
respectively by Cases $1_5$, $2_5$, $4_5$, $5_5$, $7_5$, $10_5$ and
$11_5$.\end{proof}


\section{Proof of Theorem 6}\label{secao prova teorema 6}

\begin{proof}[Proof of Theorem 6]The
direction function $H$ has a  root $Q=(q,0)$ where $q$ is given by
\eqref{eq p novo}.
Moreover, $H$ assumes positive values on the left of $Q$ and
negative values on the right of $Q$. Note that when $\alpha
\rightarrow -1$ so $Q \rightarrow \infty$ under the line $\{y=0 \}$
and it occurs the configurations showed in Theorem 4.

$\diamond$ \textit{Case $1_6$. $d<e$, Case $2_6$. $d=e$, Case $3_6$.
$d>e$, Case $4_6$. $d<S$, Case $5_6$. $d=S$ and Case $6_6$. $d>S$,
Case $7_6$. $d< h$, Case $8_6$. $d = h$, Case $9_6$. $h < d < i$,
Case $10_6$. $d=i$, Case $11_6$. $i < d < j$, Case $12_6$. $d = j$
and Case $13_6$. $j < d$:} Analogous to Cases $1_5$, $2_5$, $3_5$,
$4_5$, $5_5$, $6_5$, $7_5$, $8_5$, $9_5$, $10_5$, $11_5$, $12_5$ and
$13_5$ respectively, except that here the $\Sigma-$saddle $Q$ takes
place on the right of $d$, $e$, $S$ and $i$ when these points
appear. 

%

The bifurcation diagram is illustrated in Figure \ref{fig diagrama
bif teo 4 5 e 6}. Each topological behavior can be represented
respectively by Cases $1_6$, $2_6$, $4_6$, $5_6$, $7_6$, $10_6$ and
$11_6$.\end{proof}


\section{Proof of Theorem B}\label{secao prova teorema B}

\begin{proof}[Proof of Theorem B] Since in Equation \eqref{eq fold-sela com
parametros novos} we can take $\alpha$ in the interval
$(-1-\varepsilon_0,-1+\varepsilon_0)$ we conclude that Theorems 4, 5
and 6 prove that this equation, with $\tau=v$, unfolds generically
the (Visible) Fold$-$Saddle singularity. Its bifurcation diagram
contains all typical configurations and all distinct topological
behavior described in Theorems 4, 5 and 6. So, the number of typical
configurations is $39$ and the number of distinct topological
behavior is $21$. Moreover, each topological behavior can be
represented respectively by the Cases $1_4$,  $1_5$, $1_6$, $2_4$,
$2_5$, $2_6$,  $4_4$, $4_5$, $4_6$, $5_4$, $5_5$, $5_6$, $7_4$,
$7_5$, $7_6$, $10_4$, $10_5$, $10_6$, $11_4$, $11_5$ and $11_6$.

\begin{figure}[h!]
\begin{center}
\psfrag{A}{$1_5$}\psfrag{B}{$1_4$}\psfrag{C}{$1_6$}
\psfrag{D}{$2_6$}\psfrag{E}{$2_5$}
\psfrag{F}{$2_4$}\psfrag{G}{$3_5$}\psfrag{H}{$3_4$}
\psfrag{I}{$3_6$}\psfrag{J}{$4_6$}
\psfrag{K}{$5_6$}\psfrag{L}{$4_4$}\psfrag{M}{$4_5$}\psfrag{N}{$6_4$}
\psfrag{O}{$6_6$}
\psfrag{P}{$6_5$}\psfrag{Q}{$7_5$}\psfrag{R}{$7_4$}\psfrag{S}{$7_6$}
\psfrag{T}{$13_6$}
\psfrag{U}{$8_4$}\psfrag{V}{$12_4$}\psfrag{X}{$9_5$}\psfrag{Y}{$8_4$}
\psfrag{W}{$8_5$}
\psfrag{Z}{$9_4$}\psfrag{0}{$9_5$}\psfrag{1}{$8_6$}\psfrag{2}{$13_4$}
\psfrag{3}{$13_5$}
\psfrag{4}{$9_4$}\psfrag{5}{$9_6$}\psfrag{6}{$10_5$}\psfrag{7}{$10_4$}
\psfrag{8}{$13_6$} \psfrag{9}{$11_5$}\psfrag{`}{$\vspace{.10cm}-1<
\alpha$}\psfrag{~}{$\alpha =
-1$}\psfrag{:}{$17_5$}\psfrag{@}{$18_5$}
\psfrag{#}{$14_4$}\psfrag{^}{$17_6$}\psfrag{&}{$17_4$}\psfrag{*}{$18_6$}
\psfrag{}{$$}
\psfrag{-}{$10_6$}\psfrag{_}{$16_4$}\psfrag{=}{$19_5$}\psfrag{+}{$19_6$}
\psfrag{}{$$} \psfrag{|}{$\lambda
>0$}\psfrag{;}{$15_4$}\psfrag{,}{$\hspace{.4cm}\alpha<-1$}\psfrag{<}{$15_5$}
\psfrag{.}{$11_6$} \psfrag{>}{$11_4$}\psfrag{?}{$\beta <
0$}\psfrag{a}{$\lambda =
0$}\psfrag{b}{$\beta=0$}\psfrag{c}{$\beta>0$} \psfrag{d}{$\lambda <
0$}\psfrag{e}{$12_5$}\psfrag{f}{$12_6$}\psfrag{g}{$14_5$}\psfrag{}{$$}
\epsfxsize=13cm \epsfbox{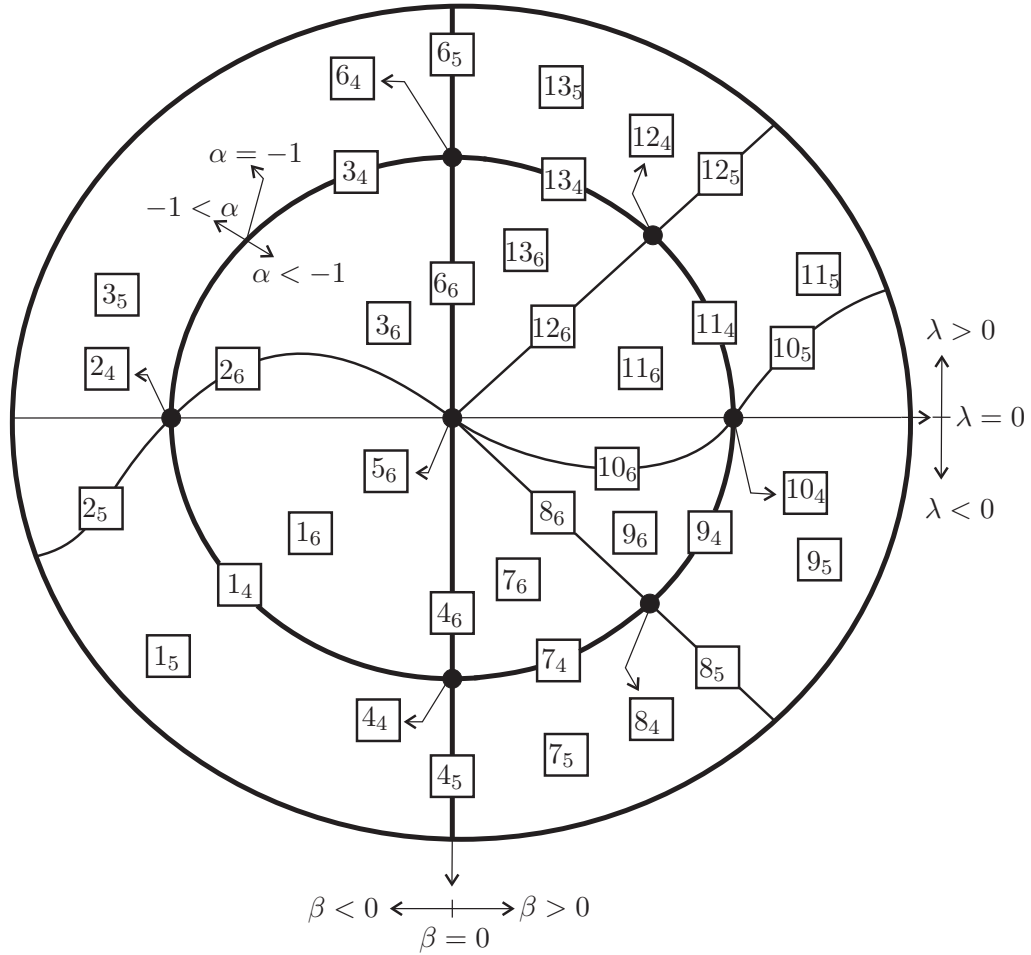}
\caption{\small{Bifurcation diagram of the (Visible) Fold$-$Saddle
singularity.}} \label{fig diagrama bifurcacao teo B}
\end{center}
\end{figure}

The full behavior of the three$-$parameter family of non$-$smooth
vector fields presenting the normal form \eqref{eq fold-sela com
parametros novos}, with $\tau=v$, is illustrated in Figure \ref{fig
diagrama bifurcacao teo B} where we consider a sphere around the
point $(\lambda , \mu, \beta) = (0,0,0)$ with a small ray and so we
make a stereographic projection defined on the entire sphere, except
the south pole. Still in relation with this figure, the numbers
pictured correspond to the occurrence of the cases described in the
previous theorems. As expected, the cases $5_4$ and $5_5$ are not
represented in this figure because they are, respectively, the
center and the south pole of the sphere.\end{proof}

\vspace{1cm}

\noindent {\textbf{Acknowledgments.} The first and the third authors
are partially supported by a FAPESP-BRAZIL grant 2007/06896-5. The
second author is partially supported by a FAPESP-BRAZIL grant
2007/08707-5. }

\end{document}